\documentclass[a4, 11pt.]{amsart}
\usepackage{amssymb,amscd, hyperref, bm, color }
\usepackage{amsmath}
\usepackage{enumerate}
\usepackage{stmaryrd}  
\usepackage[mathcal]{eucal}
\usepackage{array,float}
\usepackage{longtable}
\usepackage{multirow}
\usepackage{cite}
\usepackage[left=3cm, right=3cm]{geometry}

\usepackage{xy}
\input xy
\xyoption{all}
\usepackage{pdflscape} 
\usepackage{hyperref}
\hypersetup{
	colorlinks=blue,
	linkcolor=blue,
	filecolor=magenta,      
	urlcolor=cyan,
}

\usepackage[draft]{graphicx}

\numberwithin{equation}{section}
\newtheorem{thm}{Theorem}[section]

\newtheorem{lemma}[thm]{Lemma}

\theoremstyle{definition}

\newtheorem*{remark*}{Remark}
\newtheorem{remark}[thm]{Remark}

\newcommand{\SL}{\mathrm{SL}}
\newcommand{\GL}{\mathrm{GL}}

\newcommand{\cO}{\mathfrak{o}}
\newcommand{\mat}[4]{\left[ \begin{matrix}
		#1 & #2 \\  #3 & #4 \end{matrix}\right] }

\newcommand{\val}{\mathrm{val}}

\renewcommand{\wp}{\mathfrak{p}}
\renewcommand{\det}{{\mathrm{det}}}

\newcommand{\Char}{\mathrm{Char}}

\newcommand{\dvr}{\mathrm{DVR}}
\newcommand{\dvrtwo}{\dvr_2}
\newcommand{\dvrp}{\dvr_p}
\newcommand{\dvrtwoplus}{{\dvrtwo^+}}
\newcommand{\dvrtwozero}{{\dvrtwo^\circ}}
\newcommand{\ee}{{\mathrm{e}}}

\title[Branching rules]{Branching rules for the restriction
of regular representations  of $\GL_2(\cO/\wp^r)$ to $\SL_2(\cO/\wp^r).$
}

\author{M Hassain}
\address{ Harish-Chandra Research Institute, a CI of Homi Bhabha National Institute, Chhatnag Road, Jhunsi,
	Prayagraj - 211019, India.}

\email{hassainm@hri.res.in}



\begin{document}
	\begin{abstract}
		Let $\cO$ be a compact discrete valuation ring with maximal ideal $\wp$ such that the finite residue field $\cO/\wp$ has characteristic $p.$ 
		For $r\geq2$ and  $p=2,$ we obtain the branching rules for the restriction of a regular representation  of $\GL_2(\cO/\wp^r)$ to $\SL_2(\cO/\wp^r).$  
		These results have different behaviour than that of the known case of $p\neq2.$  

	\end{abstract}
	\maketitle
	
	
\section{Introduction}\label{sec:intro}
	For a prime $p$, let $\dvrp$ denote the set of all compact discrete valuation rings with finite residue field $\mathbb F_q$  of characteristic $p.$
	 For $\cO \in \dvrp,$ let $\wp$ be the unique maximal ideal of $\cO.$  For $r\geq 1,$ the finite quotient $\cO/\wp^r$ is denoted by $\cO_r.$
	 Let  $\dvrp^{\circ} = \{ \cO \in \dvrp \mid \Char(\cO) = 0   \}$ and  $\dvrp^{+} = \{ \cO \in \dvrp \mid \Char(\cO) = p   \}$.
	 For $\cO \in \dvrp^{\circ}$, we say $\cO$ has ramification index $\ee$ if $p \cO = \wp^\ee $. 
	Let $\GL_n(\cO_r) $ be the group of $n \times n$ invertible matrices with entries from $\cO_r$ and $\SL_n(\cO_r)$ be the subgroup of $\GL_n(\cO_r)$ consisting of all determinant one matrices.

	Recently, Patel-Singla \cite{MR4399251} proved that the restriction of a regular 
	representation $\rho$ of $\GL_n(\cO_r)$  to $\SL_n(\cO_r)$ is multiplicity free. 
	They also showed that for $(p,2)=(p,n)=1,$ $ \mathrm{Res}^{\GL_n(\cO_r)}_{\SL_n(\cO_r)}(\rho)$ has at most $n$ many irreducible constituents. In particular, for $p\neq 2,$ the representation  $ \mathrm{Res}^{\GL_2(\cO_r)}_{\SL_2(\cO_r)}(\rho)$ for a regular  representation $\rho$ of $\GL_2(\cO_r)$ has at most two irreducible constituents.
	In this article we study the restriction problem for 
	the case $p=n=2.$ In particular, we show the following result.
	For $r \geq 2,$ we denote $\ell = \lceil r/2 \rceil,$  $\ell' = \lfloor r/2 \rfloor$ and 
	$$ \mathbf N_r(\cO)=\begin{cases}
		q^{\lfloor \frac{\ell'}{2} \rfloor}& \mathrm{if} \, \, \cO \in \dvrtwoplus,\\
		q^{\lfloor \frac{\ell'}{2} \rfloor}& \mathrm{if} \, \, \cO \in \dvrtwozero \,\, \mathrm{and}\, \, \ell'< 2 \mathrm{e},\\ 
		2q^{\mathrm{e}}& \mathrm{if} \, \, \cO \in \dvrtwozero \,\, \mathrm{and}\, \, \ell'\geq 2 \mathrm{e},
	\end{cases}$$
	where $\mathrm{e}$ is the ramification index of $\cO$ if  $\cO \in \dvrtwozero.$
	\begin{thm}\label{thm:min number}
		Let $\cO \in \dvrtwo$ and $r\geq 2.$ Then there exists a regular representation $\rho$ of $\GL_2(\cO_r)$ such that the restriction $ \mathrm{Res}^{\GL_2(\cO_r)}_{\SL_2(\cO_r)}(\rho)$ has at least $\mathbf N_r(\cO)$ many irreducible constituents.
	\end{thm}
	
	In \cite{mypaper1, mypaper2}, we  constructed all 
	regular  representations of  $\SL_2( \cO_r)$ for all  $\cO \in \dvrtwoplus$ with $r \geq 1$, and for all $\cO \in \dvrtwozero$ with ramification index $\ee$ and $r \geq 2\ee.$
	By using the construction, we  proved that for $\cO \in \dvrtwozero$ with ramification index $\ee$ and $\cO' \in \dvrtwoplus,$   contrary
	to the expectation, the group algebras  $\mathbb C[\SL_2(\cO_{r})]$ and  $\mathbb C[\SL_2( \cO'_{r})]$
	are not isomorphic for any $r \geq 2 \ee +2$.
	In \cite{mypaper1}, we also studied the representation growth of $\SL_2( \cO)$ and  proved that for $\cO \in \dvrtwoplus,$ the abscissa of convergence of the representation zeta function of  $\SL_2( \cO)$ is $1,$ resolving the last remaining open
	case of this problem.

	Let $r \geq 2$ and $\pi$ be a fixed uniformizer of the ring $\cO.$ 
For $i\in\{\ell,\ell'\},$ let $M^{i}$ and $K^{i}$ be the $i^{th}$ congruence subgroups of $\GL_2(\cO_r)$ and $\SL_2(\cO_r)$ respectively. i.e., $M^{i}=I + \pi^{i} M_2(\cO_r)$ and $K^{i}=M^{i}\cap \SL_2(\cO_r).$ 
	Recall from \cite{mypaper1} that  $M^{\ell}$ is an abelian normal subgroup of $\GL_2(\cO_r)$ and the set of all one-dimensional representations of $M^{\ell}$ is given by $\widehat{M^{\ell}}=\{\psi_A \mid A \in M_2(\cO_{\ell'}) \},$ where 
	$\psi: \cO_r  \rightarrow \mathbb C^\times$ is a fixed additive  one-dimensional representation such that $\psi(\pi^{r-1}) \neq 1,$  and  $\psi_A: 
	M^\ell \rightarrow \mathbb C^\times$ is defined  by $\psi_A(I + \pi^{\ell} B) = \psi(\pi^{\ell}\mathrm{trace}(\tilde{A}B))$ for all $I + \pi^{\ell} B \in M^{\ell}.$ Here  $\tilde{A} \in M_2(\cO_r)$ is a lift of $A.$ 
	For  an irreducible representation $\rho$ of $\GL_2(\cO_r),$ it is easy to observe that   $\rho$ is regular if and only if $\langle \rho|_{M^\ell}, \psi_{A} \rangle \neq 0 $ implies $A$ is a cyclic matrix.
	Recall that a matrix $A \in M_2(\cO_r)$ for $r \geq 1,$ is called cyclic if there exists a vector $v \in \cO_r\oplus \cO_r$  such that $\{v, Av\}$ generate $\cO_r\oplus \cO_r$ as a  free $\cO_r$-module.

We obtain the following branching rules 
for the restriction problem for 
the case $p=n=2.$ For $A \in M_2(\cO_{\ell'}),$ let $C_{\GL_2(\cO_{\ell'})}(A)=\{X\in \GL_2(\cO_{\ell'}) \mid AX=XA \}.$ For a group $G$ and an irreducible representation $\phi$ of a subgroup $H$ of $G,$ let $\mathrm{Irr}(G \mid \phi)$ denote the set of all inequivalent irreducible constituents of the induced representation $\mathrm{Ind}_H^G(\phi).$

	\begin{thm}\label{thm:char zero}
		Let $\cO \in \dvrtwozero$ with   ramification index $ \mathrm{e},$   and $r\geq 4\mathrm{e}+2.$ Let $A\in M_2(\cO_{\ell'})$ be cyclic.
		Then for   $\rho\in \mathrm{Irr}(\GL_2(\cO_{r})\mid  \psi_{A}),$ 
		the following hold.
		\begin{enumerate}
			\item If $\mathrm{trace}(A)\in \cO_{\ell'}^\times,$ then  $ \mathrm{Res}_{\SL_2(\cO_{r})}^{ \GL_2(\cO_{r})}(\rho)$ is irreducible.
			\item If $\mathrm{trace}(A)\in \pi\cO_{\ell'},$ then $ \mathrm{Res}_{\SL_2(\cO_{r})}^{ \GL_2(\cO_{r})}(\rho)$ 
			is a direct sum of $\frac{(q-1)q^{\ell'-1}}{|\det(C_{\GL_2(\cO_{\ell'})}(A))|}$ many irreducible representations of dimension
			$\frac{\dim(\rho)\times|\det(C_{\GL_2(\cO_{\ell'})}(A))|}{(q-1)q^{\ell'-1}}.$
			 
		\end{enumerate}
		
	\end{thm}

	\begin{thm}\label{thm:char two inv}
		Let $\cO\in  \dvrtwoplus$ and $r\geq 2.$ Let $A\in M_2(\cO_{\ell'})$ be cyclic such that $\mathrm{trace}(A)\in \cO_{\ell'}^\times.$ Then for   $\rho\in \mathrm{Irr}(\GL_2(\cO_{r})\mid  \psi_{A}),$  the following hold.
		\begin{enumerate}
			\item If $r$ is odd, then  $ \mathrm{Res}_{\SL_2(\cO_{r})}^{ \GL_2(\cO_{r})}(\rho)$ is irreducible.
			\item If $r$ is even and $\mathrm{trace}(A)$ is not a perfect square, then   $ \mathrm{Res}_{\SL_2(\cO_{r})}^{ \GL_2(\cO_{r})}(\rho)$ is irreducible.
			\item	 If $r$ is even and $\mathrm{trace}(A)$ is a perfect square,  then   $ \mathrm{Res}_{\SL_2(\cO_{r})}^{ \GL_2(\cO_{r})}(\rho)$ is either irreducible or a direct sum of two irreducible representations of dimension $\frac{\dim(\rho)}{2}.$
		\end{enumerate}
	\end{thm}
	
	\begin{thm}\label{thm:char two non-inv}
		Let $\cO\in  \dvrtwoplus$ and $r\geq 2.$  Let $A\in M_2(\cO_{\ell'})$ be cyclic such that $\mathrm{trace}(A)\in \pi\cO_{\ell'}.$ For   $\rho\in \mathrm{Irr}(\GL_2(\cO_{r})\mid  \psi_{A}),$ let $\Delta(\rho) $ be the number of  irreducible   constituents  of
		$ \mathrm{Res}_{\SL_2(\cO_{r})}^{ \GL_2(\cO_{r})}(\rho) .$ Then the following hold for all $\rho\in \mathrm{Irr}(\GL_2(\cO_{r})\mid  \psi_{A}).$
		\begin{enumerate}
			\item If $r$ is even, then  $\frac{(q-1)q^{\ell'-1}}{|\det(C_{\GL_2(\cO_{\ell'})}(A))|}\leq  \Delta(\rho)  \leq 4\times \frac{(q-1)q^{\ell'-1}}{|\det(C_{\GL_2(\cO_{\ell'})}(A))|}. $
			\item If $r$ is odd, then  $\frac{(q-1)q^{\ell'-1}}{|\det(C_{\GL_2(\cO_{\ell'})}(A))|}\leq  \Delta(\rho) \leq  q^3\times \frac{(q-1)q^{\ell'-1}}{|\det(C_{\GL_2(\cO_{\ell'})}(A))|}. $
		\end{enumerate}
		
	\end{thm}

\section{Notations}	
In this section, we recall some notations and definitions from \cite{mypaper1,mypaper2}, which are used in this article.
For a finite group $G$ and $h\in G,$  let $C_G(h) = \{g \in G \mid gh = hg   \}$ be the centralizer of $h$ in $G.$ Similarly for any representation $\phi$ of a normal subgroup $N$ of $G,$ the group $C_G(\phi) = \{ g \in G \mid \phi^g \cong \phi \}$ be the inertia group 
of $\phi.$

Let $\cO\in  \dvrtwo$ and $r\geq 2.$
Recall from Section~\ref{sec:intro} that 	
	 we denote $\ell = \lceil r/2 \rceil$ and $\ell' = \lfloor r/2 \rfloor.$ 
	Also $\pi$ is a fixed uniformizer of the ring $\cO.$ 
	For $1\leq i\leq r,$ let $M^{i}$ and $K^{i}$ be the $i^{th}$ congruence subgroups of $\GL_2(\cO_r)$ and $\SL_2(\cO_r)$ respectively. i.e., $M^{i}=I + \pi^{i} M_2(\cO_r)$ and $K^{i}=M^{i}\cap \SL_2(\cO_r).$ 
	For  $A\in M_2(\cO_{\ell'}),$ let $\psi_{[A]}=\mathrm{Res}^{M^\ell}_{K^\ell}(\psi_{A}).$ 
	For a  cyclic $A = \mat 0 {a^{-1}\alpha } a\beta \in M_2(\cO_{\ell'}),$
	fix a lift $\tilde{A} =  \mat 0 {\tilde{a}^{-1} \tilde{\alpha}}{\tilde a}{\tilde{\beta}}  \in M_2(\cO_{r})$ of $A.$ For $i \in \{ \ell, \ell' \}$, define 
	\[ \mathrm{h}_{\tilde{A}}^{i} = \{ x \in \cO_{r } \mid 2x = 0  \,\, \mathrm{mod}\, (\pi^{i} ), \, \,  x(x+\tilde \beta) = 0 \,\, \mathrm{mod}\, (\pi^{i} ) \}.
	\]
	For $x\in \cO_r,$ let $e_x = \mat1{\tilde{a}^{-1}x}01.$ Let  $\mathrm{H}_{\tilde{A}}^{i} = \{ e_x  \mid x \in \mathrm h_{\tilde{A}}^i \}  $. Then $\mathrm{H}_{\tilde{A}}^{i} $ is an abelian group for $i \in \{ \ell, \ell' \}$. 
	By \cite[Lemma~2.2]{mypaper1}, we have 	$C_{\SL_2(\cO_{r})} (\psi_{[A]})=C_{\SL_2(\cO_{r})} (\psi_{A})\mathrm{H}_{\tilde A}^{\ell'}$ and $C_{\GL_2(\cO_{r})} (\psi_{A})=C_{\GL_2(\cO_{r})} (\tilde{A})M^{\ell'}.$ Let 
	$C_S^\ell(\tilde{A}) = ( C_{\GL_2(\cO_{r})} (\tilde{A}) M^{\ell} ) \cap \SL_2(\cO_{r})$ and  
	$D_S^{\ell}(\tilde A)= ( C_{\GL_2(\cO_{r})} (\tilde{A}) M^{\ell} ) \cap K^1.$	
	Define the following subsets of $\mathrm{H}_{\tilde A}^{\ell}$ and $\mathrm{H}_{\tilde A}^{\ell'}$:
	\begin{itemize}
		\item $\mathbb E_{\tilde{A}} := \{e_\lambda   \in \mathrm{H}_{\tilde{A}}^\ell \mid \psi_{[A]}\,\, \mathrm{extends}\,\, \mathrm{to}\,\, C_S^{\ell}(\tilde{A}) \langle e_\lambda\rangle\},$
		\item $\mathbb E^\prime_{\tilde{A}} := \{e_\lambda   \in \mathrm{H}_{\tilde{A}}^{\ell'} \mid \psi_{[A]}\,\, \mathrm{extends}\,\, \mathrm{to}\,\, D_S^{\ell}(\tilde{A}) \langle e_\lambda\rangle\},$
	\end{itemize}
where $\langle e_\lambda\rangle$ denotes the group generated by $ e_\lambda.$ 
The sets  
$\mathbb E_{\tilde{A}}$ and $\mathbb E^\prime_{\tilde{A}}$ are in bijective correspondence with the sets $ E_{\tilde{A}} := \{\lambda \in \mathrm{h}_{\tilde{A}}^\ell \mid e_\lambda \in \mathbb E_{\tilde{A}} \}$ and $ E^\prime_{\tilde{A}} := \{\lambda \in \mathrm{h}_{\tilde{A}}^{\ell'} \mid e_\lambda \in \mathbb E^\prime_{\tilde{A}} \}$ respectively.

\section{Proof of the results}	
In this section, we prove Theorems~\ref{thm:min number}, \ref{thm:char zero}, \ref{thm:char two inv} and \ref{thm:char two non-inv}.
 Let $A \in M_2(\cO_{\ell'}) $ be cyclic and $\tilde{A}$ be a lift of $A$ in $M_2(\cO_{r}).$ 
 Up to conjugation
  we can  assume 
 $A=\mat 0 {a^{-1}\alpha } a\beta$ and  $\tilde{A}=\mat 0 {\tilde{a}^{-1} \tilde{\alpha}}{\tilde a}{\tilde{\beta}},$ see \cite[Section~4]{mypaper1}.
 For $ d \in \cO^\times_{r},$ define 
 $A_d=\mat{\gamma(d)}{0}{0}{1} A \mat{\gamma(d)}{0}{0}{1}^{-1},$ where $\gamma:\cO_{r}^\times \rightarrow \cO_{\ell'}^\times$  is the natural projection.
 From \cite[Lemma~2.2]{mypaper1}, we have $C_{\GL_2(\cO_{r})} (\psi_{A})=C_{\GL_2(\cO_{r})} (\tilde{A})M^{\ell'}.$ 
Therefore it is easy to observe that $C_{\GL_2(\cO_{r})} (\psi_{A_d})=\mat{d}{0}{0}{1} C_{\GL_2(\cO_{r})} (\psi_A)\mat{d}{0}{0}{1}^{-1}.$ For $\phi\in \mathrm{Irr}(C_{\GL_2(\cO_{r})} (\psi_A)\mid  \psi_{A})$ and $d\in \cO_{r}^\times,$ define  a representation $\phi^d$ of  $ C_{\SL_2(\cO_{r})} (\psi_{A_d})=C_{\GL_2(\cO_{r})} (\psi_{A_d})\cap \SL_2(\cO_{r})$ by $\phi^d (X)=\phi(\mat{d}{0}{0}{1}^{-1} X \mat{d}{0}{0}{1}).$ 
Let $\mathcal D_A$ be a fixed set of representatives of $\cO_{r}^\times/\det (C_{\GL_2(\cO_{r})} (\psi_A)).$ 

\begin{lemma}\label{lem:res-ind-GL2}
	Let $A \in M_2(\cO_{\ell'}) $ be cyclic.
	For $\rho\in \mathrm{Irr}(\GL_2(\cO_{r})\mid  \psi_{A}),$ the following hold.
	\begin{enumerate}
		\item  There exists a unique  $\phi\in \mathrm{Irr}(C_{\GL_2(\cO_{r})} (\psi_A)\mid  \psi_{A})$ such that $\rho \cong  \mathrm{Ind}_{C_{\GL_2(\cO_{r})} (\psi_A)}^{ \GL_2(\cO_{r})}(\phi).$
		\item  $ \mathrm{Res}_{\SL_2(\cO_{r})}^{ \GL_2(\cO_{r})}(\rho) \cong \oplus_{d\in \mathcal D_A} \mathrm{Ind}_{C_{\SL_2(\cO_{r})} (\psi_{A_d})}^{ \SL_2(\cO_{r})}(\phi^d).$
	\end{enumerate}
	
\end{lemma}
\begin{proof}
	Observe that (1) directly follows from Clifford theory.  For (2), observe that the set $\{\mat{d}{0}{0}{1} \mid d \in \mathcal D_A\}$ form a set of representatives for the double cosets $\SL_2(\cO_{r}) \backslash \GL_2(\cO_{r})/ C_{\GL_2(\cO_{r})} (\psi_A).$ 
	Therefore by \cite[Proposition~22]{MR0450380} and the definition of $\phi^d$,  we obtain
	$ \mathrm{Res}_{\SL_2(\cO_{r})}^{ \GL_2(\cO_{r})}(\rho) \cong \oplus_{d\in  \mathcal D_A} \mathrm{Ind}_{C_{\SL_2(\cO_{r})} (\psi_{A_d})}^{ \SL_2(\cO_{r})}(\phi^d).$
\end{proof}
\begin{remark}\label{rem:GL2 dim 1 or q}
		Since $C_{\GL_2(\cO_{r})} (\psi_{A})=C_{\GL_2(\cO_{r})} (\tilde{A})M^{\ell'},$ 
	from \cite[Lemma~3.1]{mypaper1} and \cite[Proposition~4.11]{mypaper2} we obtain that for any  $\phi \in \mathrm{Irr}(C_{\GL_2(\cO_{r})} (\psi_A)\mid  \psi_{A}),$ $\dim(\phi)=\begin{cases}
		1&\mathrm{if}\,\, r \,\, \mathrm{even,}\\ q &\mathrm{if}\,\, r \, \,\mathrm{odd.}
	\end{cases}$
	
	\end{remark}
\begin{lemma}\label{lem:D_A cardinality}
		For cyclic $A \in M_2(\cO_{\ell'}) ,$ $| \mathcal D_A|=\frac{(q-1)q^{\ell'-1}}{|\det(C_{\GL_2(\cO_{\ell'})}(A))|}.$ Further, if   $\mathrm{trace}(A)\in \cO_{\ell'}^\times,$ then $| \mathcal D_A|=1.$
\end{lemma}
\begin{proof}
	Since $C_{\GL_2(\cO_{r})} (\psi_{A})=C_{\GL_2(\cO_{r})} (\tilde{A})M^{\ell'},$ 
	 we obtain $|\det(C_{\GL_2(\cO_{r})} (\psi_{A}))|=|\det(C_{\GL_2(\cO_{\ell'})}(A))|\times q^\ell.$ Therefore  $| \mathcal D_A|=\frac{|\cO_{r}^\times|}{|\det(C_{\GL_2(\cO_{r})} (\psi_{A}))|}=\frac{(q-1)q^{\ell'-1}}{|\det(C_{\GL_2(\cO_{\ell'})}(A))|}.$ If   $\mathrm{trace}(A)\in \cO_{\ell'}^\times,$ then \cite[Lemmas~3.4 and 4.1]{mypaper1} together give  $|\det(C_{\GL_2(\cO_{\ell'})}(A))|=|\cO_{\ell'}^\times|=(q-1)q^{\ell'-1}.$ Hence  $| \mathcal D_A|=1.$
\end{proof}

\begin{lemma}\label{lem:ind-irred-cases}
	Let $A \in M_2(\cO_{\ell'}) $ be cyclic. Let
	$\phi\in \mathrm{Irr}(C_{\GL_2(\cO_{r})} (\psi_A)\mid  \psi_{A})$ and $d\in \cO_{r}^\times.$
	Then the induced representation $\mathrm{Ind}_{C_{\SL_2(\cO_{r})} (\psi_{A_d})}^{ \SL_2(\cO_{r})}(\phi^d)$ is irreducible for the following cases.
	\begin{enumerate}
		\item $\cO \in \dvrtwozero$ and  $r\geq 4\mathrm{e}+2,$ where $ \mathrm{e}$ is the  ramification index of $\cO.$
		\item $\cO \in \dvrtwoplus,$ $r$ is odd  and $\mathrm{trace}(A)\in \cO_{\ell'}^\times.$
		\item $\cO \in \dvrtwoplus,$ $r$ is  even and $\mathrm{trace}(A)\in \cO_{\ell'}^\times\setminus (\cO_{\ell'}^\times)^2.$
	\end{enumerate}
\end{lemma}
\begin{proof}
	Assume $r\geq 2$ for $\cO \in \dvrtwoplus$ and  $r\geq 4\ee+2$ for $\cO \in \dvrtwozero$ with ramification index $\ee.$  
By definition of $\phi^d,$  to show the result, it is enough to show for $d=1.$  Note that $\phi^1$ is the restriction of $\phi$ to $C_{\SL_2(\cO_{r})} (\psi_{A}). $ 
	For odd $r,$ by Remark~\ref{rem:GL2 dim 1 or q},
	we get   $\dim(\phi^1)=q.$ Therefore by 
	\cite[Theorem~2.2]{mypaper2}, we obtain that  
	$\mathrm{Ind}_{C_{\SL_2(\cO_{r})} (\psi_{A})}^{ \SL_2(\cO_{r})}(\phi^1)$ is irreducible.
	
	For even $r,$ by Remark~\ref{rem:GL2 dim 1 or q}, we have $\dim(\phi^1)=1.$ 
	  Assume 
	  $A=\mat 0 {a^{-1}\alpha } a\beta$ and  $\tilde{A}=\mat 0 {\tilde{a}^{-1} \tilde{\alpha}}{\tilde a}{\tilde{\beta}}.$
	 To show $\mathrm{Ind}_{C_{\SL_2(\cO_{r})} (\psi_{A})}^{ \SL_2(\cO_{r})}(\phi^1)$ is irreducible, by \cite[Theorem~2.4]{mypaper1},  it is enough to show that $  C_{\SL_2(\cO_{r})} (\psi_{A})\mathbb{ E}_{\tilde A}=C_{\SL_2(\cO_{r})}(\psi_{A}).$
 	For  $\cO \in \dvrtwozero$ and  even $r\geq 4\mathrm{e}+2,$
	 by \cite[Theorem~5.6]{mypaper1}, we have 
	  $E_{\tilde{A}}=   \pi^\ell\cO_{r}.$
	   Similarly, for $\cO \in \dvrtwoplus$ and even  $r\geq 2$ such that $\beta=\mathrm{trace}(A)\in \cO_{\ell'}^\times\setminus (\cO_{\ell'}^\times)^2,$ from \cite[Proposition~3.8(1)]{mypaper1}, we obtain that $\mathrm{h}_{\tilde{A}}^\ell=\{0, \tilde{\beta}\}+\pi^\ell \cO_{r}.$ Hence by using \cite[Theorem~5.18]{mypaper1}, we have  $E_{\tilde{A}} =  \pi^\ell\cO_{r}.$  Therefore, in both cases, $  C_{\SL_2(\cO_{r})} (\psi_{A})\mathbb{ E}_{\tilde A}=C_{\SL_2(\cO_{r})}(\psi_{A}).$ Hence the result follows.
	 
	 \end{proof}

\begin{proof}[{\bf Proof of Theorem~\ref{thm:min number}}]
Let $A=\mat{0}{0}{1}{0}\in M_2(\cO_{\ell'}).$  We show that for any $ \rho \in  \mathrm{Irr}(\GL_2(\cO_{r})\mid  \psi_{A}), $ the restriction
$ \mathrm{Res}_{\SL_2(\cO_{r})}^{ \GL_2(\cO_{r})}(\rho) $ has at least $\mathbf N_r(\cO)$ many irreducible constituents. By Lemma~\ref{lem:res-ind-GL2}, it is enough to show that $|\mathcal D_A|= \mathbf N_r(\cO).$  
Since $A$ is cyclic, we have  $C_{\GL_2(\cO_{\ell'})}(A)=\{xI + yA \mid x,y \in \cO_{\ell'} \}\cap \GL_2(\cO_{\ell'}).$ Therefore by Lemma~\ref{lem:D_A cardinality}, $| \mathcal D_A|=\frac{(q-1)q^{\ell'-1}}{|\{x^2\mid x \in \cO_{\ell'}^\times \}|}=\frac{|\cO_{\ell'}^\times|}{|\{x^2\mid x \in \cO_{\ell'}^\times \}|}.$
Since the map $x\mapsto x^2$ is an endomorphism of the multiplicative group $\cO_{\ell'}^\times,$ we obtain $| \mathcal D_A|=|\{x \in \cO_{\ell'}^\times \mid x^2=1  \}|.$ For $ \cO \in \dvrtwoplus,$ since $2=0,$ it is easy to see that $\{x \in \cO_{\ell'}^\times \mid x^2=1  \}=1+\pi^{\lceil \frac{\ell'}{2} \rceil}\cO_{\ell'}.$ Hence $|\{x \in \cO_{\ell'}^\times \mid x^2=1  \}|=q^{\lfloor \frac{\ell'}{2} \rfloor}.$

For  $ \cO \in \dvrtwozero$ with ramification index $\ee,$ let $w\in \cO_{\ell'}^\times$ such that $2=\pi^\ee w.$ Note that if $x^2=1,$ then $x\in 1+\pi\cO_{\ell'}.$ For $x=1+\pi y,$ $x^2=1+2\pi y+\pi^2y^2=1+\pi^2 y (\pi^{\ee -1}w+y).$ Therefore $x^2=1$ if and only if $y (\pi^{\ee -1}w+y) \in \pi^{\ell'-2}\cO_{\ell'},$ which is equivalent to 
$\val(y)+\val(\pi^{\ee -1}w+y)\geq \ell'-2.$ Note that 
$$\{y \in \cO_{\ell'} \mid \val(y)+\val(\pi^{\ee -1}w+y)\geq \ell'-2 \}=\begin{cases}
	\{ 0, \pi^{\ee -1}w\} + \pi^{(\ell'-2)-(\ee -1)}\cO_{\ell'} & \mathrm{if}\,\, \ee -1<\lceil\frac{\ell'-2}{2}  \rceil,\\
	\pi^{\lceil\frac{\ell'-2}{2} \rceil}\cO_{\ell'} & \mathrm{if}\,\, \ee -1 \geq \lceil\frac{\ell'-2}{2} \rceil.
\end{cases} $$
It is easy to see that the conditions $\ee -1<\lceil\frac{\ell'-2}{2}  \rceil $ and $\ee -1\geq \lceil\frac{\ell'-2}{2}  \rceil $ are equivalent to  $2\ee <\ell'$ and  $2\ee \geq \ell'$ respectively. Therefore
$$\{x \in \cO_{\ell'}^\times \mid x^2=1  \}=\begin{cases}
	\{ 1, 1+\pi^{\ee}w\} + \pi^{\ell'-\ee }\cO_{\ell'} & \mathrm{if}\,\, 2\ee <\ell',\\
	1+\pi^{\lceil\frac{\ell'}{2} \rceil}\cO_{\ell'} & \mathrm{if}\,\, 2\ee \geq \ell'.
	
\end{cases} $$
Hence $|\{x \in \cO_{\ell'}^\times \mid x^2=1  \}|=\begin{cases}
		2q^{\mathrm{e}}& \mathrm{if} \, \,  2 \mathrm{e}<\ell',\\
	q^{\lfloor \frac{\ell'}{2} \rfloor}& \mathrm{if} \, \,  2 \mathrm{e}\geq \ell'.
\end{cases}$ This completes the proof of Theorem~\ref{thm:min number}.

\end{proof}

\begin{proof}[{\bf Proof of Theorem~\ref{thm:char zero}}]
	Note that for  $\phi\in \mathrm{Irr}(C_{\GL_2(\cO_{r})} (\psi_A)\mid  \psi_{A})$ and $d\in \cO_{r}^\times,$ we have $ \dim(\phi^d)=\dim(\phi)$ and $|C_{\SL_2(\cO_{r})} (\psi_{A_d})|=|C_{\SL_2(\cO_{r})} (\psi_{A})|.$ Therefore dimension of $  \mathrm{Ind}_{C_{\SL_2(\cO_{r})} (\psi_{A_d})}^{ \SL_2(\cO_{r})}(\phi^d)$ is equal to $ \frac{\dim(\phi)\times|\SL_2(\cO_{r})|}{|C_{\SL_2(\cO_{r})} (\psi_{A})|},$ which is independent of $d.$ 
	Now Theorem~\ref{thm:char zero} directly follows from Lemmas~\ref{lem:res-ind-GL2}, \ref{lem:D_A cardinality} and \ref{lem:ind-irred-cases}.
	\end{proof}

	\begin{proof}[{\bf Proof of Theorem~\ref{thm:char two inv}}]
		Note that Theorem~\ref{thm:char two inv}(1)-(2) directly follows from Lemmas~\ref{lem:res-ind-GL2}, \ref{lem:D_A cardinality} and \ref{lem:ind-irred-cases}.
		To prove Theorem~\ref{thm:char two inv}(3),  assume $r$ is even and 
		$\mathrm{trace}(A)$ is a perfect square in  $ \cO_{\ell'}^\times.$
		 By Lemma~\ref{lem:D_A cardinality},  $|\mathcal D_A|=1.$  Therefore by Lemma~\ref{lem:res-ind-GL2}, there exists $\phi \in  \mathrm{Irr}(C_{\GL_2(\cO_{r})} (\psi_A)\mid  \psi_{A})$ such that  $ \mathrm{Res}_{\SL_2(\cO_{r})}^{ \GL_2(\cO_{r})}(\rho) \cong  \mathrm{Ind}_{C_{\SL_2(\cO_{r})} (\psi_{A})}^{ \SL_2(\cO_{r})}(\phi^1).$ By Remark~\ref{rem:GL2 dim 1 or q}, 
		 $\phi$ is one dimensional. Hence $\phi^1\in \mathrm{Irr}(C_{\SL_2(\cO_{r})} (\psi_A)\mid  \psi_{[A]}).$ 
		 Since $\mathrm{trace}(A) \in \cO_{\ell'}^\times,$ by using \cite[Lemma~2.2 and Proposition~3.8(1)]{mypaper1}, it is easy to see that $[C_{\SL_2(\cO_{r})} (\psi_{[A]}):C_{\SL_2(\cO_{r})} (\psi_{A})]=2.$ Therefore, by Clifford theory,  the induced representation $\mathrm{Ind}_{C_{\SL_2(\cO_{r})} (\psi_{A})}^{ \SL_2(\cO_{r})}(\phi^1)$ is either irreducible or a direct sum of two irreducible representations of same dimension. Hence  Theorem~\ref{thm:char two inv}(3) holds.
	\end{proof}
	
	To prove Theorem~\ref{thm:char two non-inv}, we need the following lemma.
\begin{lemma}\label{lem:min dim trace non-inv}
		Let $\cO\in  \dvrtwoplus$ and $r> 2$ be odd.  Let $A\in M_2(\cO_{\ell'})$ be cyclic such that $\mathrm{trace}(A)\in \pi\cO_{\ell'}.$ For   any $\chi\in \mathrm{Irr}(\SL_2(\cO_{r})\mid  \psi_{[A]}),$  $\dim(\chi) \geq \frac{1}{q^2}\times \frac{|\SL_2(\cO_{r})|}{|C_{\SL_2(\cO_{r})} (\psi_{A})|} .$
\end{lemma}	
\begin{proof}
	Assume 
	$A=\mat 0 {a^{-1}\alpha } a\beta$ and  $\tilde{A}=\mat 0 {\tilde{a}^{-1} \tilde{\alpha}}{\tilde a}{\tilde{\beta}}.$
	 Let $\chi\in \mathrm{Irr}(\SL_2(\cO_{r})\mid  \psi_{[A]}).$ 
	 By 
	\cite[Theorem~2.3]{mypaper2}, there exist  
	a subgroup $H_\chi$ of $ C_{\SL_2(\cO_{r})} (\psi_{[A]})$ containing 
	$D_S^{\ell}(\tilde A)= ( C_{\GL_2(\cO_{r})} (\tilde{A}) M^{\ell} ) \cap K^1,$
	and 
	an extension $\phi$ of $\psi_{[A]}$ to $H_\chi$ such that $\chi \cong \mathrm{Ind}_{H_\chi}^{ \SL_2(\cO_{r})}(\phi).$ Therefore, to show the result, it is enough to show that $|H_\chi|\leq q^2 \times |C_{\SL_2(\cO_{r})} (\psi_{A})|. $   Note that by \cite[Lemma~2.2]{mypaper1}, 
	$C_{\SL_2(\cO_{r})} (\psi_{[A]})=C_{\SL_2(\cO_{r})} (\psi_{A})\mathrm{H}_{\tilde A}^{\ell'}.$
	Let $m_\chi=\{\lambda\in \mathrm{h}_{\tilde{A}}^{\ell'} \mid Xe_\lambda \in H_\chi \, \, \mathrm{for} \, \, \mathrm{some}\, \, X\in C_{\SL_2(\cO_{r})} (\psi_{A})\}.$ Since $C_{\SL_2(\cO_{r})} (\psi_{A}) \trianglelefteq C_{\SL_2(\cO_{r})} (\psi_{A})\mathrm{H}_{\tilde A}^{\ell'}$
	and $\{e_\lambda \mid \lambda \in \pi^{\ell'}\cO_r \} \subseteq C_{\SL_2(\cO_{r})} (\psi_{A}),$
	it is easy to observe that $m_\chi$ form an additive subgroup of $ \mathrm{h}_{\tilde{A}}^{\ell'}$ such that $\pi^{\ell'}\cO_r \subseteq m_\chi.$  This implies
	\begin{equation}\label{eqn:asw1}
		|H_\chi|\leq |C_{\SL_2(\cO_{r})} (\psi_{A})\{e_\lambda \mid \lambda \in  m_\chi \}  | \leq  |C_{\SL_2(\cO_{r})} (\psi_{A})|\times \frac{|m_\chi |}{|\pi^{\ell'}\cO_r|}.
	\end{equation}

	We claim that $ m_\chi \subseteq E_{\tilde A}'.$ Let $\lambda \in m_\chi.$
	 Choose $X\in C_{\SL_2(\cO_{r})} (\psi_{A})$ 
	 such that  $Xe_\lambda \in H_\chi.$ Since $C_{\SL_2(\cO_{r})} (\psi_{A})$ stabilizes $\mathrm{Res}_{D_S^\ell(\tilde{A})}^{H_\chi}(\phi)$ (from \cite[Lemma~4.4]{mypaper2}), we obtain that $e_\lambda=X^{-1}(Xe_\lambda)$ stabilizes $\mathrm{Res}_{D_S^\ell(\tilde{A})}^{H_\chi}(\phi).$  Note that the quotient  group $D_S^\ell(\tilde{A})\langle e_\lambda\rangle/D_S^\ell(\tilde{A})$ is cyclic.  Therefore by \cite[Corollary~11.22]{MR2270898},  
	 $\mathrm{Res}_{D_S^\ell(\tilde{A})}^{H_\chi}(\phi)$ extends to $D_S^\ell(\tilde A) \langle e_\lambda  \rangle  .$  In particular
	 $\psi_{[A]}$ has an extension to $D_S^\ell(\tilde A) \langle e_\lambda  \rangle  .$ 
	  Therefore  $ \lambda \in E_{\tilde A}'.$ Hence the claim.

	 Consider the following set of valuations:
	 $$\val(E_{\tilde A}' \setminus \pi^{\ell'}\cO_r)=\{\val(\lambda) \mid \lambda\in E_{\tilde A}'\setminus \pi^{\ell'}\cO_r \}.$$
	 From \cite[Theorem~3.10, Lemma~3.11]{mypaper2} and \cite[Proposition~5.16(1)-(2)]{mypaper1}, we obtain $|\val(E_{\tilde A}' \setminus \pi^{\ell'}\cO_r)|\leq 2.$ Let $\val(E_{\tilde A}' \setminus \pi^{\ell'}\cO_r)\subseteq \{v_1,v_2\}$ for some $v_1<v_2<\ell'.$ Then we must have 
	 \begin{equation}\label{eqn:asd2}
	 	E_{\tilde A}' \subseteq \pi^{v_1 }\cO_r^\times \cup \pi^{v_2 }\cO_r^\times \cup\pi^{\ell'}\cO_r .
	 \end{equation}

	 Let $m_\chi^\prime = m_\chi \cap \pi^{v_2 }\cO_r.$ Then $\pi^{\ell'}\cO_r \subseteq m_\chi^\prime \subseteq m_\chi \subseteq E_{\tilde A}'.$ We use that $m_\chi^\prime$ is a group and (\ref{eqn:asd2}) to observe that if any two elements of $m_\chi^\prime$ are equal modulo $\pi^{v_2 +1 }\cO_r$ then these must be equal modulo $\pi^{\ell'}\cO_r.$ Therefore $[m_\chi^\prime :\pi^{\ell'}\cO_r ]\leq q.$  
	 By similar argument we get $[m_\chi :m_\chi^\prime ]\leq q.$  Therefore
	 $$ \frac{|m_\chi |}{|\pi^{\ell'}\cO_r|}\leq [m_\chi :m_\chi^\prime ]\times [m_\chi^\prime :\pi^{\ell'}\cO_r ] \leq q^2 . $$
	 This with (\ref{eqn:asw1}) give $|H_\chi|\leq q^2 \times |C_{\SL_2(\cO_{r})} (\psi_{A})|. $ 
	 
\end{proof}

\begin{proof}[{\bf Proof of Theorem~\ref{thm:char two non-inv}}]
		Note that  Lemmas~\ref{lem:res-ind-GL2} and \ref{lem:D_A cardinality} together give  $\frac{(q-1)q^{\ell'-1}}{|\det(C_{\GL_2(\cO_{\ell'})}(A))|}\leq  \Delta(\rho)$ for any $r\geq 2.$ So we focus on upper bound of $\Delta(\rho).$ By Lemma~\ref{lem:res-ind-GL2}, there exists 
		 $\phi\in \mathrm{Irr}(C_{\GL_2(\cO_{r})} (\psi_A)\mid  \psi_{A})$ such that $\rho \cong  \mathrm{Ind}_{C_{\GL_2(\cO_{r})} (\psi_A)}^{ \GL_2(\cO_{r})}(\phi)$ 
		and  $ \mathrm{Res}_{\SL_2(\cO_{r})}^{ \GL_2(\cO_{r})}(\rho) \cong \oplus_{d\in \mathcal D_A} \mathrm{Ind}_{C_{\SL_2(\cO_{r})} (\psi_{A_d})}^{ \SL_2(\cO_{r})}(\phi^d).$ If $r$ is even, then  
		 by Remark~\ref{rem:GL2 dim 1 or q}, 	$\phi$ is one dimensional. Therefore for each $d \in \mathcal D_A,$ the representation $\phi^d \in \mathrm{Irr}(C_{\SL_2(\cO_{r})} (\psi_{A_d})\mid  \psi_{[A_d]}).$ 
		   This together with \cite[Theorems~2.4(3) and 5.20(1)]{mypaper1} imply that  the induced representation 
		$\mathrm{Ind}_{C_{\SL_2(\cO_{r})} (\psi_{A_d})}^{ \SL_2(\cO_{r})}(\phi^d_s)$ decomposes into maximum four irreducible components. Therefore $\Delta(\rho) \leq 4\times |\mathcal D_A| = 4\times \frac{(q-1)q^{\ell'-1}}{|\det(C_{\GL_2(\cO_{\ell'})}(A))|}, $ where the last equality follows from Lemma~\ref{lem:D_A cardinality}.
		
		For odd $r,$ by Remark~\ref{rem:GL2 dim 1 or q}, we have $\dim(\phi)=q.$  Therefore for each $d \in \mathcal D_A,$ the induced representation  $\mathrm{Ind}_{C_{\SL_2(\cO_{r})} (\psi_{A_d})}^{ \SL_2(\cO_{r})}(\phi^d)$ has dimension $q\times \frac{|\SL_2(\cO_{r})|}{|C_{\SL_2(\cO_{r})} (\psi_{A_d})|}.$
		Since $A_d$ is cyclic and $\mathrm{trace}(A_d)= \mathrm{trace}(A)\in \pi\cO_{\ell'},$ by Lemma~\ref{lem:min dim trace non-inv}, every $\chi \in \mathrm{Irr}(\SL_2(\cO_{r}) \mid  \psi_{[A_d]})$ has dimension greater than or equal to $\frac{1}{q^2}\times \frac{|\SL_2(\cO_{r})|}{|C_{\SL_2(\cO_{r})} (\psi_{A_d})|}.$
		 Therefore  the induced representation 
		$\mathrm{Ind}_{C_{\SL_2(\cO_{r})} (\psi_{A_d})}^{ \SL_2(\cO_{r})}(\phi^d_s)$ decomposes into maximum $q^3$ many irreducible components. Thus $\Delta(\rho)  \leq q^3\times |\mathcal D_A|= q^3 \times \frac{(q-1)q^{\ell'-1}}{|\det(C_{\GL_2(\cO_{\ell'})}(A))|}. $ 
		Hence Theorem~\ref{thm:char two non-inv} holds.
	
	\end{proof}

		\section*{Acknowledgements} 
	
	The author immensely thank Pooja Singla for several helpful discussions and comments on this work.

	\bibliography{refs}{}
	\bibliographystyle{siam}
	
\end{document}